\documentclass[12pt]{amsart}

\pdfoutput=1

\usepackage[UKenglish]{babel}
\usepackage{amssymb}
\usepackage{amsmath}
\usepackage{graphicx}
 \usepackage{hyperref}

\newtheorem{theorem}{Theorem}[section]
\newtheorem{proposition}[theorem]{Proposition}

\newtheorem{lemma}[theorem]{Lemma}

\theoremstyle{definition}
\newtheorem{example}[theorem]{Example}
\theoremstyle{remark}
\newtheorem*{remark}{Remark}

\def\E{\mathbb{E}}
\def\PP{\mathbb{P}}

\def\wP{\widetilde{\mathbb{P}}}

\def\RR{\mathbb{R}}
\def\NN{\mathbb{N}}
\def\loc{L^1_\textrm{loc}}
\def\cF{\mathcal{F}}

\def\d{\text{\rm d}}

\def\1{ 1\hspace{-2.9pt}{\rm l}}

\begin{document}


\title{Martingale property of generalized stochastic exponentials}

\author{Aleksandar Mijatovi\'c}
\address{Department of Mathematics, Imperial College London }
\email{a.mijatovic@imperial.ac.uk}

\author{Nika Novak}
\address{Faculty of Mathematics and Physics, University of Ljubljana
\newline \indent and Department of Mathematics, Imperial College London}
\email{nika.novak@fmf.uni-lj.si}

\author{Mikhail Urusov}
\address{Institute of Mathematical Finance, Ulm University, Germany}
\email{mikhail.urusov@uni-ulm.de}

\subjclass[2000]{60G44, 60G48, 60H10, 60J60.}
\keywords{Generalized stochastic exponentials; local martingales vs. true martingales; one-dimensional diffusions.}

\begin{abstract}
For a real Borel measurable function
$b$,
which satisfies certain integrability conditions,
it is possible to define a stochastic integral
of the process
$b(Y)$
with respect to a Brownian motion 
$W$,
where 
$Y$
is a diffusion driven by
$W$.
It is well know that the stochastic exponential of this 
stochastic integral is a local martingale.
In this paper we consider the case of 
an arbitrary Borel measurable function
$b$
where it may not be possible to define the stochastic 
integral of 
$b(Y)$
directly.
However the notion of the stochastic exponential can be generalized.
We define a non-negative process
$Z$,
called \textit{generalized stochastic exponential},
which is not necessarily a local martingale.
Our main result gives deterministic necessary and sufficient conditions 
for 
$Z$
to be a local, true or uniformly integrable martingale. 
\end{abstract}

\maketitle

\section{Introduction}
A \emph{stochastic exponential} of 
$X$
is a process $\mathcal{E}(X)$ defined by 
$$
\mathcal{E}(X)_t = \exp\left\{X_t - X_0 - \frac 1 2 \langle X \rangle_t\right\}
$$
for some continuous local martingale $X$, 
where 
$\langle X \rangle$ 
denotes a quadratic variation of $X$.
It is well known that the process 
$\mathcal{E}(X)$ 
is also a continuous local martingale.
The characterisation of the martingale property of 
$\mathcal{E}(X)$
has been studied extensively in the literature
because this question appears naturally in many situations.

In the case of one dimensional processes, necessary and sufficient conditions for the process $\mathcal{E}(X)$ 
to be a martingale were recently studied by Engelbert and Senf in \cite{ES}, Blei and Engelbert in \cite{BE} and 
Mijatovi\'c and Urusov in \cite{MU09}. 
In \cite{ES} $X$ is a general continuous local martingale and the characterisation is given in terms of the Dambis-Dubins-Schwartz time-change 
that turns 
$X$ 
into a Brownian motion. 
In \cite{BE} $X$ is a strong Markov continuous local 
martingale and the condition is deterministic, 
expressed in terms of the speed measure of $X$.

In \cite{MU09} the local martingale 
$X$ 
is of the form 
$X_t= \int_0^t b\left(Y_u\right) \d W_u$ 
for some measurable function 
$b$
and a one-dimensional diffusion 
$Y$ 
with drift 
$\mu$ 
and volatility 
$\sigma$ 
driven by a Brownian motion 
$W$. 
In order to define the stochastic integral 
$X$,
an assumption that the function 
$\frac{b^2}{\sigma^2}$ 
is locally integrable on the entire state space of the process 
$Y$ 
is required.
Under this restriction  the characterization of the martingale property
of 
$\mathcal{E}(X)$
is studied in~\cite{MU09},
where the necessary and sufficient 
conditions 
are deterministic and are expressed in terms of functions 
$\mu, \sigma$ and $b$   
only.

In the present paper we consider an arbitrary Borel measurable function $b$. 
In this case the stochastic integral 
$X$ 
can only be defined on some subset of 
the probability space. 
However, it is possible to 
define a non-negative 
possibly discontinuous process 
$Z$,
known as a generalized stochastic exponential,
on the entire probability space. 
It is a consequence of the definition that, 
if the function
$b$
satisfies the required local integrability condition, 
the process
$Z$
coincides with
$\mathcal{E}(X)$.
We show that 
the process 
$Z$ 
is not necessarily a local martingale. 
In fact 
$Z$
is a local martingale if and only if it is continuous.
We find a deterministic necessary and sufficient condition for 
$Z$ 
to be a local martingale,
which is expressed in terms of 
local integrability of the quotient
$\frac{b^2}{\sigma^2}$ multiplied by a linear function. 
We also characterize the processes 
$Z$ 
that are true martingales 
and/or 
uniformly integrable martingales. 
All the necessary and sufficient conditions are
deterministic and are given in terms of functions 
$\mu, \sigma$ and $b$.

The paper is structured as follows.  
In Section 2 we define the notion of generalized stochastic exponential
and study its basic properties. 
The main results are stated in Section 3, where we give a 
necessary and sufficient condition for the process $Z$ defined 
by \eqref{Z} and \eqref{Z_inf} to be a local martingale,  
a true martingale or a uniformly  integrable martingale. 
Finally, in Section 4 we prove Theorem~\ref{deterministic} 
that is central in obtaining  
the deterministic characterisation of the martingale property
of the process $Z$. 
Appendix A contains an auxiliary fact that is used in Section~2.

\section{Definition of Generalized Stochastic Exponential}
Let $J=(l,r)$ be our state space, where $-\infty \le l < r \le \infty$. Let us define a $J$-valued diffusion $Y$ on a probability space $(\Omega,\cF, (\cF_t)_{t \in [0, \infty)}, \PP)$ driven by a stochastic differential equation
\begin{equation*}\label{main1}
\d Y_t = \mu\left(Y_t\right) \d t + \sigma \left(Y_t\right) \d W_t, \quad Y_0=x_0 \in J,
\end{equation*}
where $W$ is  a $(\cF_t)$-Brownian motion and  $\mu$ and $\sigma$ real, Borel measurable functions defined on $J$ that satisfy the Engelbert-Schmidt conditions
\begin{align}\label{m2}
 \sigma (x) \ne 0 \quad \forall x \in J,\\ \label{m3}
 \frac{1}{\sigma^2}, \frac{\mu}{\sigma^2} \in \loc (J).
\end{align}
With $\loc(J)$ we denote the class of locally integrable functions, i.e. real Borel measurable functions defined on $J$ that are integrable on every compact subset of $J$. Engelbert-Schmidt conditions guarantee existence of a weak solution that might exit the interval $J$ and is unique in law (see \cite[Chapter 5]{KS}). Denote by $\zeta$ the exit time of $Y$. In addition, we assume that the boundary points are absorbing, i.e. the solution $Y$ stays at the boundary point at which it exits on the set $\{\zeta < \infty\}.$
Let us note that we assume that $(\cF_t)$ is generated neither by $Y$ nor by $W$.

We would like to define a process $X$ as a stochastic integral of a process $b(Y)$ with respect to Brownian motion $W$, where $b:J \to \RR$ is an arbitrary Borel measurable function.   Before further discussion, we should establish if the stochastic integral can be defined.

Define a set
$$
A = \{ x \in J ;\  \textstyle\frac{b^2}{\sigma^2} \not\in \loc (x)\},
$$ 
where $\loc (x)$ denote a set of real, Borel measurable functions $f$ such that $\int_{x-\varepsilon}^{x+\varepsilon}f(y) \d y < \infty$ for some $\varepsilon > 0$. Then $A$ is closed and its complement is a union of open intervals. Let $\tau_A= \inf\{t \ge 0;\  Y_t \in A \}$ and $\zeta^A= \zeta \wedge \tau_A$. 
Then
$$
\int_0^t b^2\left(Y_u\right) \d u < \infty \ \ \PP \text{-a.s. on } \{t < \zeta^A\}.
$$
This follows from Proposition \ref{det1} and the fact that a continuous process $Y$ on $\{t < \zeta^A\}$ reaches only values in  an open interval that is a component of the complement of $A$, where $\frac{b^2}{\sigma^2}$ is locally integrable.  

Let us define
$
A_{n} = \{ x \in J; \ \rho (x, A \cup \{l,r\}) \le \frac 1 n \},
$ 
where $\rho(x,y) = |\arctan x - \arctan y|, x,y \in \bar J$, and
set $\zeta_n^A= \inf \{t \ge 0; \ Y_t \in A_n \}$. Since $\zeta_n^A < \zeta^A$ on the set $\{\zeta^A < \infty\}$, we have $\int_0^{t \wedge \zeta_n^A} b^2(Y_u) \d u < \infty \ \PP$-a.s. Thus, we can define the stochastic integral $\int_0^{t \wedge \zeta_n^A} b(Y_u) \d W_u$ for every $n$. Since the integrals $\int_0^{t \wedge \zeta_n^A} b(Y_u) \d W_u$ and $\int_0^{t \wedge \zeta_{n+1}^A}b(Y_u) \d W_u$ coincide on $\{t < \zeta_{n}^A\}$ and $\zeta_n^A \uparrow \zeta^A$, we can define $\int_0^{t \wedge \zeta^A} b(Y_u) \d W_u$ as a limit of the integrals $\int_0^{t \wedge \zeta_n^A} b(Y_u) \d W_u$.

In the case where $A$ is not empty or $Y$ exits the interval $J$, the stochastic exponential cannot be defined. However, we can define a \emph{generalized stochastic exponential}  $Z$ in the following way for every $t \in [0, \infty)$
\begin{align}\label{Z}
Z_t &= \left\{   
\begin{array}{ll}
 \exp\{\int_0^{t} b(Y_u) \d W_u - \frac 1 2 \int_0^{t} b^2 (Y_u) \d u \} &, t < \zeta^A \\
 \exp\{\int_0^{\zeta} b(Y_u) \d W_u - \frac 1 2 \int_0^{\zeta} b^2 (Y_u) \d u \} &, t \ge \zeta^A = \zeta, \int_0^{\zeta} b^2\left(Y_u\right) \d u < \infty\\
 0 &, t \ge \zeta^A = \tau_A \text{ or } \\
   &\  t \ge \zeta^A= \zeta, \int_0^{\zeta} b^2 \left(Y_u\right) \d u = \infty 
\end{array}
\right.
\end{align}
The different behaviour of $Z$ on $\{t \ge \zeta^A = \zeta\}$ and $\{t \ge \zeta^A = \tau_A\}$ follows from the fact, that after the exit time $\zeta$ the process $Y$ is stopped, while  this does not happen after $\tau_A$. From the definition of the set $A$ the integral $\int_0^t b^2\left(Y_u\right) \d u$ is infinite for every $t > \tau_A$. Therefore, we set $Z=0$ on the set $\{t \ge \zeta^A = \tau_A\}$.
 
Let us define the processes 
\begin{equation}\label{Z_bar}
\bar Z_t = \exp\left\{\int_0^{t \wedge\zeta^A}b\left(Y_u\right) \d W_u - \frac 1 2 \int_0^{t \wedge \zeta^A} b^2\left(Y_u \right) \d u \right\},
\end{equation}
where we set $\bar Z_t =0$ for $t \ge \zeta^A$ on $\{\zeta^A < \infty, \int_0^{\zeta^A} b^2(Y_u) \d u = \infty\},$
and
$$
S_t = \exp\left\{ \int_0^{\tau_A} b\left(Y_u\right) \d W_u - \frac 1 2 \int_0^{\tau_A} b^2\left(Y_u\right)\d u\right\} \1_{\{t \ge \zeta^A = \tau_A, \int_0^{\tau_A} b^2\left(Y_u\right)\d u < \infty\}}.
$$
Then we can write 
\begin{equation}\label{Z_eq}
Z= \bar Z - S.
\end{equation}
Now $Z$ is not necessarily a continuous process.
Furthermore,  $\bar Z$ is positive local martingale and therefore a supermartingale. The process $S$  has increasing paths. Hence,
$$
\E[Z_t | \cF_s] \le \bar Z_s - \E[S_t | \cF_s] \le \bar Z_s - S_s=Z_s.
$$
It follows that $Z$ is a supermartingale and we can define
\begin{align}
 \label{Z_inf}
Z_{\infty} = \lim_{t \to \infty} Z_t.
\end{align}

\begin{remark}
Note that we should not use (\ref{Z}) for $t=\infty$ because in (\ref{Z}) $Z_{\infty}$ is not well defined  on $\{\zeta^A = \infty\}$.
\end{remark}

We may assume that $x_0 \not\in A$. Otherwise, $Z \equiv 0$ and hence it is a martingale. 

A path of the process $Z$ defined by (\ref{Z}) and (\ref{Z_inf}) is equal to a path of a stochastic exponential if $\zeta^A = \infty$. Otherwise, if $\zeta^A < \infty$, it has one of the following forms:
\begin{enumerate}
 \item[$(i)$]
 $\tau_A < \zeta$ and $\int_0^{\tau_{A}} b^2(Y_t)\d t < \infty$ (see Figure \ref{slikca1});
   \begin{figure}[h]

   \includegraphics[width=3.3in]{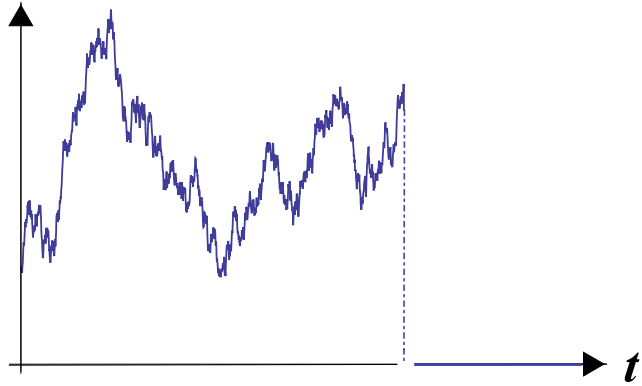}
    \caption{If $\tau_A<\zeta$, then the process $Z$ is positive up to time $\tau_A$ and is equal to zero afterwards. If the integral $\int_0^{\tau_A} b^2(Y_t)\d t$ is finite, then $Z_t$ approaches a positive value as $t$ approaches $\tau_A$. Therefore,  there is a jump at $t= \tau_A$.}\label{slikca1}
   \end{figure}

 \item[$(ii)$]
 $\zeta^A < \infty$ and $\int_0^{\zeta^{A}} b^2(Y_t) \d t = \infty$ (see Figure \ref{slikca2});
   \begin{figure}[h]

    \includegraphics[width=3.3in]{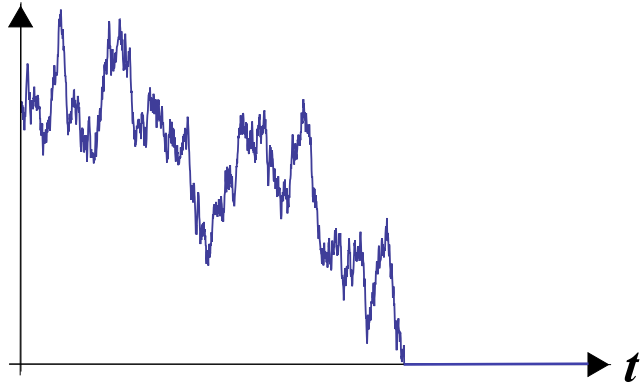}
     \caption{If $\zeta^A < \infty$ and $\int_0^{\zeta^A} b^2(Y_t) \d t = \infty$, then the process $Z$ is zero after the time $\zeta^A$. Since the limit of $Z_t$ is zero as $t$ approaches $\zeta^A$, there is no jump.}\label{slikca2}
   \end{figure}

 \item[$(iii)$]
 $\zeta  < \tau_A$ and $\int_0^{\zeta} b^2 (Y_t) \d t < \infty$ (see Figure \ref{slikca3}).
   \begin{figure}[h]

    \includegraphics[width=3.3in]{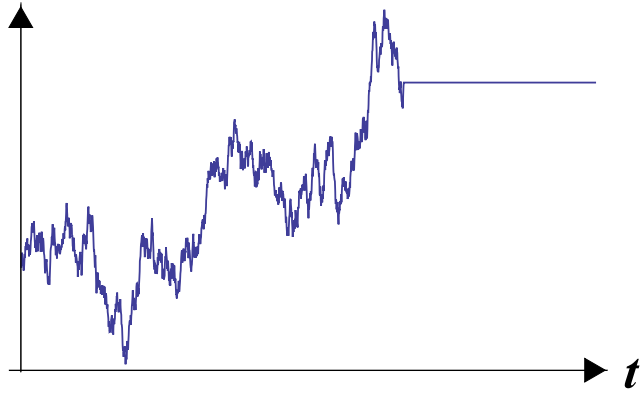}
    \caption{If $\zeta < \tau_A$, the process $Z$ is stopped after the exit time. Since $\int_0^{\zeta} b^2 (Y_t) \d t$ is finite, $Z_t$ is equal to a positive constant for $t \ge \zeta$. }\label{slikca3}
   \end{figure}

\end{enumerate}

\section{Main Results}
The case $A = \emptyset$ was studied by Mijatovi\'c and Urusov in \cite{MU09}. We generalize their result for the case where $A \ne \emptyset$. 

\subsection{The Case $A= \emptyset$}\label{prvi}
In this case we have 
\begin{equation} \label{m1}
\frac{b^2}{\sigma^2} \in \loc (J). 
\end{equation}
 The generalized stochastic exponential $Z$ defined by (\ref{Z}) and (\ref{Z_inf}) can now be written  as
$$
Z_t = \exp\left\{\int_0^{t \wedge \zeta} b\left(Y_u\right)\d W_u - \frac 1 2 \int_0^{t \wedge \zeta} b^2\left(Y_u\right) \d u\right\},
$$
where we set $Z_t=0$ for $t \ge \zeta$ on $\{\zeta < \infty, \int_0^{\zeta} b^2\left(Y_u\right)\d u = \infty\}$.
Note that in this case 
$Z$
is a local martingale.

Let us now define an auxiliary $J$-valued diffusion $\widetilde Y$ governed by the SDE
\begin{equation*}\label{main2}
\d \widetilde Y_t = \left(\mu + b \sigma\right) \left(\widetilde Y_t\right) \d t + \sigma\left(\widetilde Y_t\right) \d \widetilde W_t, \quad \widetilde Y_0 = x_0,
\end{equation*}
on some probability space $(\widetilde \Omega, \widetilde \cF, (\widetilde \cF_t)_{t \in [0, \infty)}, \wP)$. The coefficients $\mu + b \sigma$ and $\sigma$ satisfy Engelbert-Schmidt conditions since $\frac{b}{\sigma} \in \loc (J)$ (this follows from (\ref{m1})). Hence the SDE has a unique in law, possibly explosive, weak solution. As with diffusion $Y$, we denote by $\widetilde
\zeta$ the exit time of $\widetilde Y$ and assume that the boundary points are absorbing.

For an arbitrary $c \in J$ we define the scale functions $s, \widetilde s$ and their derivatives $\rho, \widetilde \rho$:
\begin{equation}\label{split}
\begin{split}
 &\rho (x) = \exp \left\{- \int_c^x \frac{2 \mu (y)}{\sigma^2(y)} \; \d y \right\}, \quad x \in  J,\\
&\widetilde \rho (x) = \rho (x) \exp \left\{ - \int_c^x \frac{2 b(y)}{\sigma(y)} \; \d y \right\}, \quad  x\in J, \\
&s(x) = \int_c^x \rho \left(y\right) \d y, \quad x \in \bar J,\\ 
&\widetilde s (x) = \int_c^x \widetilde \rho \left(y\right) \d y, \quad x \in \bar J.
\end{split}
\end{equation}
Denote by $\loc(r-)$ the set of all Borel measurable functions $f:J \to \RR$ such that $\int_{r-\varepsilon}^r \left|f\left(x\right)\right| \d x$ is finite for some $\varepsilon >0$. Similarly, we denote by $\loc(l+)$ the set of all Borel measurable functions such that $\int_l^{l+\varepsilon} \left|f\left(x\right)\right| \d x$ is finite for some $\varepsilon >0$.

We say that the endpoint $r$ is \emph{good} if
\begin{equation*}\label{good1}
s(r) < \infty \quad \text{ and }\quad \frac{(s(r) - s)b^2}{\rho \sigma^2} \in \loc (r-). 
\end{equation*}
It is equivalent to show that
\begin{equation*}\label{good3}
\widetilde s(r) < \infty \quad \text{ and }\quad \frac{(\widetilde s(r) - \widetilde s)b^2}{\widetilde \rho \sigma^2} \in \loc (r-).
\end{equation*}  
The endpoint $l$ is \emph{good} if
\begin{equation*}\label{good2}
s(l) > - \infty \quad \text{ and } \quad \frac{(s - s(l))b^2}{\rho \sigma^2} \in \loc (l+),
\end{equation*}
or equivalently
\begin{equation*}\label{good4}
\widetilde s(l) > - \infty \quad \text{ and } \quad \frac{(\widetilde s -\widetilde s(l))b^2}{\widetilde\rho \sigma^2} \in \loc (l+).
\end{equation*}
If an endpoint is not good, we say it is \emph{bad}. The good and bad endpoints were introduced in \cite{MU09}, where one can also find the proof of equivalences above.

\bigskip

We will use the following terminology:
\begin{center} 
$\widetilde Y$ \emph{exits at $r$} means $\wP (\widetilde\zeta < \infty, \lim_{t \uparrow \widetilde\zeta} \widetilde Y_t=r)>0$;

$\widetilde Y$ \emph{exits at} $l$ means $\wP (\widetilde\zeta < \infty, \lim_{t \uparrow \widetilde\zeta} \widetilde Y_t=l)>0$.
\end{center}

Define 
\begin{align} \label{v1}
&\widetilde v(x) = \int_c^x \frac{\widetilde s(x) - \widetilde s(y)}{\widetilde \rho(y)\sigma^2(y)} \; \d y, \quad x \in J, \\
\intertext{and}
&\widetilde v(r) = \lim_{x \uparrow r} \widetilde v(x), \quad \widetilde v(l) = \lim_{x \downarrow l} \widetilde v(x). \label{v2}
\end{align}

\emph{Feller's test for explosions} (see \cite[Chapter 5, Theorem 5.29]{KS}) tells us that:  
\begin{enumerate}
\item[$(i)$]
$\widetilde Y$ exits at the boundary point $r$ if and only if 
$$
\widetilde v(r) < \infty.
$$
It is equivalent to check (see \cite[Chapter 4.1]{CE})
$$
\widetilde s(r) < \infty\quad \text{ and } \quad\frac{\widetilde s(r) - \widetilde s}{\widetilde\rho \sigma^2} \in \loc (r-);
$$
\item[$(ii)$]
$\widetilde Y$ exits at the boundary point $l$ if and only if
$$
\widetilde v(l) > - \infty,
$$
which is equivalent to 
$$
\widetilde s(l) > - \infty\quad \text{ and }\quad \frac{\widetilde s - \widetilde s(l)}{\widetilde \rho \sigma^2} \in \loc (l+).
$$ 
\end{enumerate}

\begin{remark}
The endpoint $r$ (resp. $l$) is bad whenever one of the processes $Y$ and $\widetilde Y$ exits at $r$ (resp. $l$) and the other does not. 
\end{remark}

\begin{theorem}\label{MU1}
Let the functions $\mu, \sigma$ and $b$ satisfy conditions $(\ref{m2})$, $(\ref{m3})$ and $(\ref{m1})$. Then the process $Z$ is a martingale if and only if $\widetilde Y$ does not exit at the bad endpoints.
\end{theorem}

\begin{theorem}\label{MU2}
Let the functions $\mu, \sigma$ and $b$ satisfy conditions $(\ref{m2})$, $(\ref{m3})$ and $(\ref{m1})$. Then $Z$ is a uniformly integrable martingale if and only if one of the conditions $(a) - (d)$ below is satisfied:
\begin{enumerate}
 \item[$(a)$] $b=0$ a.e. on $J$ with respect to the Lebesgue measure;
 \item[$(b)$] $r$ is good and $\widetilde s (l) = -\infty$;
 \item[$(c)$] $l$ is good and $\widetilde s (r) = \infty$;
 \item[$(d)$] $l$ and $r$ are good.
\end{enumerate} 
\end{theorem}

\subsection{The Case $A \ne \emptyset$}
The following example shows that even when $A$ is not empty we can get a martingale or a uniformly integrable martingale  defined by (\ref{Z}) and (\ref{Z_inf}).
\begin{example}
 $(i)$ 
 Let us consider the case $J=\RR$, $\mu=0, \sigma=1$ and $b(x)=\frac{1}{x}$. Then $A=\{0\}$ and $Y_t =W_t, W_0=x_0 > 0$. Using It\^o's formula and the fact that Brownian motion does not exit at infinity, we get for $t<\tau_0$
\begin{align*}
Z_t 
&= \exp\left\{\int_0^t \frac{1}{W_u}\;\d W_u - \frac 1 2 \int_0^t \frac{1}{W_u^2}\; \d u\right\}\\
&=\frac{1}{x_0}W_{t}
\end{align*}
and $Z_t=0$ for $t \ge \tau_0$. Hence, $Z_t= \frac{1}{x_0}W_{t \wedge \tau_0}$
 that is a martingale.
 
 $(ii)$
Using the same functions $\mu, \sigma$ and $b$ as above on a state space $J= (- \infty, x_0+1)$ we get
$$
Z_t = \frac{1}{x_0} W_{t \wedge \tau_{0, x_0+1}},
$$
which is a uniformly integrable martingale.
\end{example}

Define maps $\alpha$ and $\beta$ on $J \setminus A$ so that 
\begin{equation}\label{endpoint}
\alpha (x), \beta (x) \in A \cup \{l,r\}\quad \text{ and } \quad x\in (\alpha(x), \beta(x)) \subset J \setminus A.
\end{equation}
So, $\alpha (x)$ is the point in $A$ that is closest to $x$ from the left side and $\beta(x)$ is the closest point in $A$ from the right side. 
Then $\frac{b^2}{\sigma^2} \in \loc (\alpha(x), \beta(x))$.
Therefore, on $(\alpha(x), \beta(x))$ functions $\mu, \sigma$ and $b$ satisfy the same conditions as in previous subsection.

We can define an auxiliary diffusion $\widetilde Y$ with values in $(\alpha(x_0), \beta(x_0))$ driven by the SDE
$$
\d \widetilde Y_t = \left(\mu + b \sigma\right)\left(\widetilde Y_t\right) \d t + \sigma\left(\widetilde Y_t\right) \d \widetilde W_t, \quad \widetilde Y_0 = x_0,
$$ 
on some probability space $(\widetilde{\Omega}, \widetilde\cF, (\widetilde\cF_t)_{t \in [0, \infty)}, \widetilde\PP)$.
There exists a unique weak solution of this equation since coefficients satisfy the Engelbert-Schmidt conditions.

As in the previous subsection we can define good and bad endpoints. We say that the endpoint $\beta(x_0)$ is \emph{good} if
\begin{equation*}
s(\beta(x_0)) < \infty \quad \text{ and }\quad \frac{(s(\beta(x_0)) - s)b^2}{\rho \sigma^2} \in \loc (\beta(x_0)-). 
\end{equation*}
It is equivalent to show and sometimes easier to check that
\begin{equation*}
\widetilde s(\beta(x_0)) < \infty \quad \text{ and }\quad \frac{(\widetilde s(\beta(x_0)) - \widetilde s)b^2}{\widetilde \rho \sigma^2} \in \loc (\beta(x_0)-).
\end{equation*}  
The endpoint $\alpha(x_0)$ is \emph{good} if
\begin{equation*}
s(\alpha(x_0)) > - \infty \quad \text{ and } \quad \frac{(s - s(\alpha(x_0)))b^2}{\rho \sigma^2} \in \loc (\alpha(x_0)+),
\end{equation*}
or equivalently
\begin{equation*}
\widetilde s(\alpha(x_0)) > - \infty \quad \text{ and } \quad \frac{(\widetilde s -\widetilde s(\alpha(x_0)))b^2}{\widetilde\rho \sigma^2} \in \loc (\alpha(x_0)+).
\end{equation*}
If an endpoint is not good, we say it is \emph{bad}.

\begin{remark}
Functions $\rho, \widetilde \rho, s, \widetilde s$ and $\widetilde v$ are defined as in (\ref{split}), (\ref{v1}) and (\ref{v2}). Now, we only need to take $c$ from the interval $(\alpha(x_0), \beta(x_0))$.
\end{remark}

Define a set
$$
B=\left\{x \in A; \int_0^{\tau_x}b^2\left(Y_t\right) \d t = \infty \ \ \PP\text{-a.s.}\right\},
$$ 
where $\tau_x = \inf \{t \ge 0; Y_t =x\}$.
The following theorem characterizes the set $B$ in a deterministic way.
\begin{theorem}\label{deterministic}
Let $\alpha$ be the function defined by (\ref{endpoint}), $\alpha(x_0) > l$  and let us write shortly $\alpha = \alpha(x_0)$. Then:
 \begin{enumerate}
  \item[$(a)$] $(x- \alpha) \frac{b^2}{\sigma^2} (x) \in \loc (\alpha + ) \; \Longleftrightarrow \; \int_0^{\tau_{\alpha}} b^2\left(Y_t \right) \d t < \infty \; \PP$-a.s. on $\{\tau_{\alpha}= \tau_A < \infty\}$;
  \item[$(b)$] $ (x- \alpha) \frac{b^2}{\sigma^2} (x) \notin \loc (\alpha + ) \; \Longleftrightarrow \; \int_0^{\tau_{\alpha}} b^2\left(Y_t\right) \d t = \infty \; \PP$-a.s. on $\{\tau_{\alpha}= \tau_A < \infty\}$.
 \end{enumerate}
\end{theorem}

Note that the assertions $(a)$ and $(b)$ in Theorem \ref{deterministic} are not the negation of each other. If the integral $\int_0^{\tau_{\alpha}} b^2\left(Y_t\right) \d t$ is not finite $\PP$-a.s. on $\{\tau_{\alpha}= \tau_A < \infty \}$, then it is infinite on some subset of $\{\tau_{\alpha}= \tau_A < \infty \}$ with positive probability. Observe that $\PP(\tau_{\alpha} = \tau_A < \infty) > 0$.

Clearly, Theorem \ref{deterministic} has its analogue for $\beta(x_0)<r$. 

Now we can show when a generalized stochastic exponential is a local martingale and when it is a true martingale.
\begin{theorem}\label{general1}
$(i)$ The generalized stochastic exponential $Z$ is a local martingale if and only if $\alpha (x_0), \beta(x_0) \in B \cup \{l,r\}$.

$(ii)$ The generalized stochastic exponential $Z$ is a martingale if and only if $Z$ is a local martingale and 
at least one of the conditions (a)-(b) below is satisfied and at least one of the conditions (c)-(d) below is satisfied:
\begin{enumerate}
 \item[(a)] $\widetilde Y$ does not exit at $\beta (x_0)$, i.e. $\widetilde v(\beta (x_0)) = \infty$ or equivalently, 
 $$
 \widetilde s (\beta (x_0)) = \infty \quad\text{ or } \quad \left(\widetilde s (\beta (x_0)) < \infty \text{ and } \frac{\widetilde s(\beta (x_0)) - \widetilde s}{\widetilde \rho \sigma^2} \notin \loc(\beta (x_0)-)\right);
$$
\item[(b)] $\beta(x_0)$ is good, 
 \item[(c)] $\widetilde Y$ does not exit at $\alpha(x_0)$, i.e. $\widetilde v(\alpha(x_0)) =- \infty$ or equivalently, 
 $$
 \widetilde s (\alpha(x_0)) = - \infty \quad\text{ or } \quad \left(\widetilde s (\alpha(x_0)) > -\infty \text{ and } \frac{\widetilde s - \widetilde s(\alpha(x_0))}{\widetilde \rho \sigma^2} \notin \loc(\alpha(x_0)+)\right);
$$
\item[(d)] $\alpha(x_0)$ is good. 
\end{enumerate}
\end{theorem}
\begin{remark}
Part $(ii)$ of Theorem \ref{general1} says that $Z$ is a martingale if and only if  
the $(\alpha(x_0), \beta(x_0))$-valued process $\widetilde Y$ can exit only at the good endpoints.
\end{remark}
\begin{proof}
$(i)$ We can write $Z=\bar Z -S$ as in (\ref{Z_eq}). Since $\left(\int_0^{t \wedge \zeta^A} b\left(Y_u\right) \d W_u\right)_t$ is a continuous local martingale, the process $\bar Z$ is a continuous local martingale. Suppose that $Z$ is a local martingale. Then $S$ can be written as a sum of two local martingales and  therefore, it is also a local martingale. It follows that $S$ is a supermartingale (since it is positive). Since $\zeta^A > 0$ and $S_0=0$, $S$ should be almost surely equal to $0$.
This happens if and only if $\alpha(x_0), \beta(x_0) \in B \cup \{l,r\}$.

$(ii)$  To get at least a local martingale  $S$ needs to be zero  $\PP$-a.s. Then $Z =\bar Z.$ Since the values of $Y$ on $[0, \zeta^A)$ do not exit the interval $(\alpha(x_0), \beta(x_0)),$ the conditions of Theorem \ref{MU1} are satisfied and the result follows.
\end{proof}

Similarly, we can characterize uniformly integrable martingale. We can use characterization in Theorem \ref{MU2} for the  process $\bar Z$ defined by (\ref{Z_bar}). As above, for $\alpha(x_0), \beta(x_0) \in B \cup \{l,r\} $  the process $Z$ defined by (\ref{Z}) and (\ref{Z_inf}) coincides with $\bar Z$. Otherwise, $Z$ is not even a local martingale.  
\begin{theorem}\label{general2}
The process $Z$ is a uniformly integrable martingale if and only if $Z$ is a local martingale and at least one of the conditions $(a)-(d)$ below is satisfied:
\begin{enumerate}
 \item[$(a)$] $b =0$ a.e. on $(\alpha(x_0), \beta(x_0))$ with respect to the Lebesgue measure;
 \item[$(b)$] $\alpha(x_0)$ is good and $\widetilde s (\beta(x_0)) = \infty$;
 \item[$(c)$] $\beta(x_0)$ is good and $\widetilde s(\alpha(x_0)) = - \infty$;
 \item[$(d)$] $\alpha(x_0)$ and $\beta(x_0)$ are good.
\end{enumerate}
\end{theorem}

\begin{remark}
If $\alpha (x_0) \in B$, then $\alpha(x_0)$ is not a good endpoint. Indeed, if $s(\alpha(x_0)) > -\infty$, then we can write
$$ 
\frac{(s(x) - s(\alpha(x_0))) b^2(x)}{\rho (x) \sigma^2(x)} = 
\frac{(s(x) - s(\alpha(x_0)))}{(x - \alpha(x_0))\rho(x)} (x - \alpha(x_0))\frac{b^2}{\sigma^2}(x).
$$
The first fraction is bounded away from zero, since it is continuous for $x > \alpha(x_0)$ and has a limit equal to $1$ as $x$ approaches $\alpha(x_0)$. Therefore, $\frac{(s - s(\alpha(x_0)))b^2}{\rho \sigma^2} \not\in \loc (\alpha(x_0) +)$.

Similarly, $\beta(x_0) \in B$ implies that $\beta(x_0)$ is not a good endpoint.

This remark simplifies the application of Theorems \ref{general1} and \ref{general2} in specific situations.
\end{remark}

\section{Proof of Theorem \ref{deterministic}}

For the proof of Theorem \ref{deterministic} we first consider the case of Brownian motion. Let $W$ be a Brownian motion with $W_0=x_0$.
Denote by $L_t^y(W)$ a local time of $W$ at time $t$ and level $y$.  
Let $-\infty < \alpha < x_0$ and consider a Borel function $b:(\alpha, \infty) \to \RR$ that is locally integrable on the interval $(\alpha, \infty)$.

\begin{lemma}\label{lemma 1}
If $(x - \alpha)b^2(x) \in \loc (\alpha +)$, then  $\int_0^{\tau_{\alpha}} b^2\left(W_t\right) \d t < \infty \; \PP$-a.s.
\end{lemma}
\begin{proof}
Let $(\beta_n)$ be an increasing sequence such that $x_0 < \beta_n < \infty$ and $\beta_n \uparrow \infty$.
By \cite[Chapter VII, Corollary 3.8]{RY} we get
\begin{align*}
 \E_{x_0} \left[\int_0^{\tau_{\alpha}\wedge \tau_{\beta_n}} b^2 \left(W_t\right)\d t\right]
 = 2 {\textstyle \frac{\beta_n-x_0}{\beta_n - \alpha}} \int_{\alpha}^{x_0} (y - \alpha) b^2\left(y\right) \d y  +
    2 {\textstyle\frac{x_0 - \alpha}{\beta_n - \alpha}} \int_{x_0}^{\beta_n} (\beta_n-y) b^2\left(y\right) \d y
\end{align*}
for every $\beta_n$.
Both integrals are finite since $b^2 \in L_{loc}^1 (\alpha, \infty)$ and $ (x - \alpha)b^2(x) \in \loc (\alpha +).$
Thus, we have $ \E_{x_0} [\int_0^{\tau_{\alpha}\wedge \tau_{\beta_n}} b^2 \left(W_t\right) \d t]< \infty$  and therefore $\int_0^{\tau_{\alpha} \wedge \tau_{\beta_n}} b^2\left(W_t\right)  \d t < \infty$ almost surely for every $n$.

Since  $\{\tau_{\alpha} < \tau_{\beta_n}\} \uparrow \{ \tau_{\alpha} < \infty\}$ almost surely as $n$ tends to infinity and 
$\PP(\{ \tau_{\alpha} < \infty\}) =1$, we get
$$
\int_0^{\tau_{\alpha}} b^2\left(W_t\right)\d t < \infty \;\; \PP\text{-a.s., } 
$$
which concludes the proof.
\end{proof}

\begin{lemma}\label{lemma 2}
 If $\int_0^{\tau_{\alpha}} b^2 \left(W_t\right) \d t < \infty$ on a set $U$ with $\PP (U)>0$, then $ (x-\alpha) b^2(x) \in \loc (\alpha +)$.
\end{lemma}
\begin{proof}
 The idea of the proof comes from \cite{J}. 
 Using the occupation times formula we can write
\begin{align*}
\int_0^{\tau_{\alpha}} b^2\left(W_t\right)\d t 
    &= \int_{\alpha}^{\infty} b^2(y) L_{\tau_{\alpha}}^y\left(W\right) \d y\\
    &\ge \int_{\alpha}^{x_0} b^2 (y) L_{\tau_{\alpha}}^y\left(W\right) \d y.
\end{align*}

Let us define a process $R_y= \frac{1}{y-\alpha}L_{\tau_{\alpha}}^y(W)$. Then  $R$ is positive and we have
\begin{equation}\label{in1}
\int_0^{\tau_{\alpha}} b^2\left(W_t\right) \d t \ge \int_{\alpha}^{x_0} R_y (y - \alpha) b^2\left(y\right) \d y.
\end{equation}
By \cite[Chapter VI, Proposition 4.6]{RY}, Laplace transform of $R_y$ is
$$
\E [\exp\{-\lambda R_y\}] = \frac{1}{1+ 2 \lambda} \quad \text{for every } y.
$$
Hence, every random variable $R_y$ has exponential distribution with $\E [R_y] =  2$.
 
Denote by $L$ an indicator function of a measurable set. We can write
$$
\E [L R_y] = \E \left[L \int_0^{\infty} \1_{\{R_y > u\}}  \d u\right] = \int_0^{\infty} \E[L \1_{\{R_y > u\}}] \d u.
$$
By Jensen's inequality we get a lower bound for the integrand 
\begin{align*}
\E [L \1_{\{R_y >u\}}]
        &= \E [(L- \1_{\{R_y \le u\}})^+]\\
        &\ge ( \E [L] - \PP [R_y \le u])^+\\
        &= (\E [L] + e^{-\frac u 2} -1)^+.
\end{align*}
Hence,
\begin{equation}\label{in2}
\E [L R_y] \ge \int_0^{\infty} (\E [L] + e^{- \frac u 2}- 1)^+\d u =C, 
\end{equation}
where $C$ is a strictly positive constant if $\E[L]$ is strictly positive.

Then we choose $L$, so that $\E [L \int_0^{\tau_{\alpha}} b^2\left(W_t\right) \d t]$ is finite. Using Fubini's Theorem and inequalities (\ref{in1}) and (\ref{in2}), we get 
\begin{align*}
 \E \left[L \int_0^{\tau_{\alpha}} b^2 \left(W_t\right) \d t\right] 
           &\ge  \int_{\alpha}^{x_0} \E [L R_y] (y- \alpha) b^2 \left(y\right) \d y\\
           &\ge C \int_{\alpha}^{x_0} (y - \alpha) b^2\left(y\right) \d y.
\end{align*}
Therefore, $(y-\alpha) b^2(y) \in \loc (\alpha +)$ if we can find an indicator function $L$ such that $\E [L]$ is strictly positive and
$\E [L \int_0^{\tau_{\alpha}} b^2\left(W_t\right) \d t]$ is finite. 

Since $\int_0^{\tau_{\alpha}} b^2\left(W_t\right) \d t < \infty$ on a set with positive measure, such $L$ exists. Indeed, denote by $L_n$ an indicator function of the set  $U_n=\{\int_0^{\tau_{\alpha}} b^2\left(W_t\right) \d t \le n\}$. Then $\E [L_n \int_0^{\tau_{\alpha}} b^2\left(W_t\right) \d t] < \infty$ for every integer $n$.  Since the sequence $(U_n)_{n \in \NN}$ is increasing, $U \subseteq \bigcup_{n \in \NN} U_n$ and $\PP(U) >0$, there exists an integer $N$ such that $\PP(U_N)>0$ and therefore $\E[L_N] > 0$.
\end{proof}

Now we return to the setting of Section~2.

\begin{proof}[Proof of Theorem \ref{deterministic}]
First, suppose that $\mu=0$ and $\sigma =1$. In this case  our diffusion $Y_t$ is equal to a (possibly stopped) Brownian motion $W_t$ with $W_0=x_0$.
The equivalences in $(a)$ and $(b)$ follow from Lemmas \ref{lemma 1} and \ref{lemma 2}.

\medskip

Now suppose that $\mu=0$ and $\sigma$ is arbitrary. 
Since $Y_t$ is a continuous local martingale, by Dambis--Dubins--Schwartz we have $Y_t = B_{\langle Y \rangle_t}$ for a Brownian motion $B$ with $B_0 = Y_0$.
Using the substitution $s = \langle Y \rangle_t$, we get
\begin{align*}
 \int_0^{\tau_{\alpha}}b^2\left(Y_t\right) \d t 
 &= \int_0^{\tau_{\alpha}} \frac{b^2}{\sigma^2} \left(Y_t\right) \d \langle Y \rangle_t\\
 &= \int_0^{\langle Y \rangle_{\tau_{\alpha}}} \frac{b^2}{\sigma^2} \left(B_s\right) \d s.
\end{align*} 
Since $B_{\langle Y \rangle_{\tau_{\alpha}}} = Y_{\tau_{\alpha}} = \alpha$ and $\langle Y \rangle_{\tau_{\alpha}} = \inf \{ s \ge 0; B_s = \alpha \}$, we can use the first part of the proof to show the assertions.

\bigskip

It only remains  to prove the general case when both $\mu$ and $\sigma$ are arbitrary. Let $Z_t = s(Y_t)$, where $s$ is the scale function of $Y$. 
Then  $Z$ satisfies SDE
$$
\d Z_t = \widetilde \sigma \left(Z_t\right) \d W_t,
$$   
where $\widetilde \sigma (x) = s'(q(x))\sigma(q(x))$ and $q$ is the inverse of $s$.

Define $\widetilde b = b \circ q$. Since $s$ is increasing and $Z_{\tau_{\alpha}} = s(Y_{\tau_{\alpha}}) = s (\alpha)$, we can also show that $\tau_{\alpha}(Y) = \tau_{s(\alpha)}(Z).$ 
Then we have
$$
\int_0^{\tau_{\alpha}(Y)} b^2\left(Y_t\right) dt = \int_0^{\tau_{s(\alpha)}(Z)} \widetilde b^2\left(Z_t\right) dt.
$$
Besides,
\begin{align*}
\int_{s(\alpha)}^{s(\alpha+ \varepsilon)} \frac{\widetilde b^2(x)}{\widetilde \sigma^2(x)} \left(x - s(\alpha)\right) \d x 
&= \int_{\alpha}^{\alpha + \varepsilon} \frac{b^2(y)}{\sigma^2(y)} (y- \alpha) \frac{(s(y)-s(\alpha))s'(y)}{ y- \alpha}\;\d y.
\end{align*} 
Fraction $\frac{(s(y)-s(\alpha))s'(y)}{ y- \alpha}$ is continuous  for $y > \alpha$ and has a positive limit  in $\alpha$. Hence it is bounded and bounded away from zero.  
It follows that $(x- \alpha)\frac{b^2}{\sigma^2}(x) \in L_{loc}^1 (\alpha+)$ if and only if $(x- s(\alpha)) \frac{\widetilde b^2}{\widetilde \sigma^2}(x) \in L_{loc}^1 (s(\alpha)+).$ 
Then the result follows from the second part of the proof.
\end{proof}

\appendix
\section{}

Let $Y$ be a $J$-valued diffusion with a drift $\mu$ and volatility $\sigma$ that satisfy Engelbert-Schmidt conditions. Let $b: J\to \RR$ be a Borel-measurable function and 
let $(c,d) \subseteq J$.

\begin{proposition}\label{det1}
A condition
$$
\frac{b^2}{\sigma^2} \in \loc (c,d) 
$$
is equivalent to
$$
\int_0^t b^2 \left(Y_u\right) \d u < \infty \; \PP\text{-a.s.  on } \{ t < \tau_{c,d} \}.
$$
\end{proposition}
\begin{proof}
Using the occupation times formula we get
\begin{align*}
 \int_0^t b^2(Y_u)  \d u 
 &= \int_0^t \frac{b^2}{\sigma^2}(Y_u) \d \langle Y \rangle_u \\
 &= \int_c^d \frac{b^2}{\sigma^2}(y) L_t^y(Y) \d y.
\end{align*}
Suppose first that $\frac{b^2}{\sigma^2} \in \loc (c,d)$. Since the function $y \mapsto L_t^y(Y)$ is c\'adl\'ag (see \cite[Chapter VI, Theorem 1.7]{RY}), it is bounded for every $t < \tau_{c,d}$ and has a compact support in $(c,d)$. The implication follows directly.

Suppose now that $\frac{b^2}{\sigma^2} \notin \loc (c,d)$. Then there exists such $\alpha \in (c,d)$ that we  have
\begin{align}\label{d1}
&\int_{\alpha- \varepsilon}^{\alpha + \varepsilon} \frac{b^2}{\sigma^2}(y)\d y = \infty \ \text{ for all } \varepsilon >0.
\end{align}
It is well known that $\PP(\tau_{\alpha} < \infty) >0$, where $\tau_{\alpha}=\{t \ge 0;\  Y_t = \alpha\}$.
By \cite[Theorem 2.7]{CE},  we have 
$$
L_t^{\alpha} (Y) >0 \text{ and } \lim_{y \uparrow \alpha} L_t^y(Y) > 0 \ \ \PP \text{-a.s.}
$$
for any $t \ge 0$ on the set $\{t > \tau_{\alpha}\}$.
Then there exists $\varepsilon >0$, such that the function $y \mapsto L_t^y(Y)$ is bounded away from zero $\PP$-a.s. on $\{t > \tau_{\alpha}\}$ on the interval
$(\alpha - \varepsilon, \alpha + \varepsilon)$. It follows that $\int_0^t b^2\left(Y_u\right) \d u = \infty $ on $\{t > \tau_{\alpha}\}$, which proves the assumption.
\end{proof}



\end{document}